**Title**

**A New Concept of optimal control for epidemic spreading by Vaccination: Technique for Assessing social optimum employing Pontryagin's Maximum Principle**

Md. Mamun-Ur-Rashid Khan [a,c], Jun Tanimoto [a,b]

[a]Interdisciplinary Graduate School of Engineering Sciences, Kyushu University, Kasuga-koen, Kasuga-shi, Fukuoka 816-8580, Japan

[b]Faculty of Engineering Sciences, Kyushu University, Kasuga-koen, Kasuga-shi, Fukuoka 816-8580, Japan

[c]Department of Mathematics, University of Dhaka, Dhaka-1000, Bangladesh

**Corresponding author**
Md. Mamun-Ur-Rashid Khan, mamun.math@du.ac.bd

**Abstract**

This research introduces a new approach utilizing optimal control theory (OCT) to assess the Social Optimum (SO) of a vaccination game, navigating the intricate considerations of cost, availability, and distribution policies. By integrating an SIRS/V epidemic model with a behavior model, the study analyzes individual vaccination strategies. A unique optimal control framework, centered on vaccination costs, is proposed, diverging significantly from previous methods. Our findings confirm the effectiveness and feasibility of this approach in managing vaccination strategies. Moreover, we examine the underlying social dilemma of the vaccination game, investigating key parameters. By calculating the Nash equilibrium (NE) through the behavior model and determining the SO using our approach, we measure the Social Efficiency Deficit (SED), quantifying the overall cost gap between the NE and SO. Results indicate that an increased waning immunity rate exacerbates the social dilemma, although higher vaccination costs partially mitigate it. This research provides valuable insights into optimizing vaccination strategies amidst complex societal dynamics.

**Highlights**

• An innovative approach to optimal control is introduced to determine the social optimum.
• A comparison between the traditional cost-based behavior model and the newly proposed optimal control model is shown.
• The proposed optimal control concept proves more reliable and effective in achieving the desired vaccination strategy.

**1. Introduction**

Vaccination is key in managing infectious diseases, yet initial shortages occur in pandemics like COVID-19. The availability of vaccines becomes vital once developed. However, during COVID-19, only wealthier nations could offer widespread vaccination, leaving low-income countries grappling with distribution challenges. Cost also impacts coverage significantly; lower costs allow authorities to provide subsidies or free distribution [1–3]. High costs present obstacles in widespread vaccine distribution, impacting individuals' decisions to participate in vaccination programs, influenced by factors like infection rates, vaccine effectiveness, and cost. Meanwhile, authorities tasked with vaccine provision aim to develop cost-efficient distribution strategies [3–6].

To implement any intervention, such as vaccination, treatment, quarantine, or isolation, in an epidemic model, it is essential to examine the model from two perspectives: the situation without the intervention and the changes that occur when the intervention is introduced [7–11]. Numerous studies have been conducted during the COVID-19 pandemic, employing various interventions in epidemic models. Many of these studies utilized an epidemic model based on the Susceptible, Infected, Recovered (SIR) framework, augmented with additional compartments such as Exposed (E), Hospitalized (H), Aware (A), Unaware (U), Treatment (T), Quarantine (Q), Protected (P), Death (D), and others [5–9,11–37]. These investigations aimed to analyze the impact of interventions on disease spread. Many studies utilized Pontryagin's Maximum Principle, an optimal control theory approach, to minimize an objective function. However, they overlooked vaccination costs, assuming the objective function only comprised infection prevalence and intervention rates squared, which doesn't reflect total social costs. To truly minimize social costs, including disease and intervention costs, a new framework must be established. This study aims to address this gap.

Furthermore, when analyzing the social dilemma within an epidemic model, it is essential to consider the costs associated with infection and vaccination [38]. Calculating the social optimum is another significant aspect of such analysis [4,38–46]. Some previous studies have employed time-constant vaccination rates to determine the social optimum, which can sometimes be impractical for authorities as it may suggest vaccinating 100% of individuals at the outset, which is nearly impossible in real-world scenarios [2–4,39,40,47–56]. Furthermore, the waning rate of immunity added more urgency to the situation. In this regard, optimal control theory provides a more suitable and mathematically acceptable approach to calculating the socially optimal vaccination level at any given time.

Our study focuses on a simplified SIRS/V epidemic model, examining how individuals respond to vaccination costs. We determine the optimal vaccination level using optimal control theory, integrating vaccination costs. Our novel objective function, aligning with Pontryagin's maximum principle, combines vaccination and infection costs. Individuals decide on vaccination based on observed infection rates and vaccination costs. Social optimum is chosen to minimize the total social cost, including infection and vaccination costs. We compare the payoff disparities between models to illustrate the Social Efficiency Deficit (SED), revealing the social dilemma within our proposed model [3,39,47,50,52].

## 2. Model Depiction

### 2.1 Epidemic model with behavior dynamics

This addresses the novelty of our concept for building the objective function for an optimal control problem; we presume a relatively simple vaccination game where a simple compartment model is coupled with a behavior model.

Our research utilized an epidemic model with four compartments based on the dynamics of SIRS/V populations. The total population is initially categorized as the susceptible group ($S$), consisting of individuals who are susceptible to the infection. These individuals can contract the infection disease determined by the transmission rate ($\beta$) and transition to the infected compartment ($I$). Subsequently, infected individuals recover from the infection at a rate of ($\gamma$) and move to the recovered compartment. To incorporate vaccination into the model, we introduced a separate compartment for vaccinated individuals ($V$). This compartment represents individuals who have received the vaccine. The individuals' transition from the susceptible compartment to the vaccination compartment is denoted as $x(t)$. This rate is determined using the behavior model [57−61], which considers factors such as the number of infected individuals in different states and the cost associated with vaccination. It is worth noting that we consider vaccination to be regarded, where we represent the vaccine's efficacy as $\eta$. Therefore, individuals in the vaccination compartment can still contract the infection at a rate of $(1-\eta)\beta$. Additionally, we accounted for the waning rate of immunity, represented by $\omega$, which captures the gradual decrease in immunity over time. The flow diagram and formulation of the proposed model are as follows:

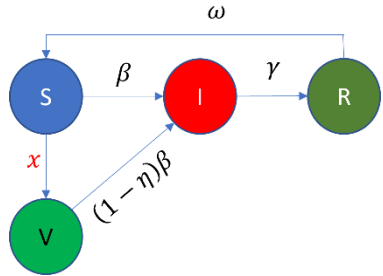

**Figure 1:** Model Flowchart (including behavioral dynamics). In this diagram, susceptible individuals may become infected at a transmission rate, $\beta$, transitioning from compartment $S$ to compartment $I$. Infected individuals can recover at a rate of $\gamma$, moving to the recovered compartment $R$. Susceptible individuals can transition to the vaccinated compartment $V$ at a rate determined by behavior dynamics, denoted by $x$. Vaccinated individuals can be infected with a discounted transmission rate $(1-\eta)\beta$ and move to the infected compartment. Recovered individuals may become susceptible again at a rate of waning immunity, represented by $\omega$.

$$\dot{S} = -\beta SI - xS + \omega R, \tag{1}$$

$$\dot{V} = xS - (1-\eta)\beta VI, \tag{2}$$

$$\dot{I} = \beta SI + (1-\eta)\beta VI - \gamma I, \tag{3}$$

$$\dot{R} = \gamma I - \omega R, \tag{4}$$

$$S(t) + V(t) + I(t) + R(t) = 1, \tag{5}$$

$$\dot{x} = mx(1-x)(cI - kc_V), \tag{6}$$

where $m$ is the inertial effect of the vaccination, $c$ is the disease cost due to infection, $c_v$ is the vaccination cost, and $k$ is the relative sensitivity due to the cost of vaccination.

### 2.2 Epidemic model with optimal control

In this section, we extend the epidemic model discussed earlier by introducing a control variable $(u)$ that represents the flux of vaccination from the susceptible compartment to the vaccination compartment. We utilize optimal control theory to find the optimal value $u$, specifically applying Pontryagin's maximum principle. This principle helps us determine the optimal control strategy that minimizes the objective function, considering the dynamics of the epidemic model and the constraints imposed by the system. By employing this approach, we can identify the most effective vaccination strategy to combat the spread of the disease and maximize the desired outcomes. All the model parameters and their description are shown in Table 1. The schematic diagram and formulation of the model can be summarized as follows:

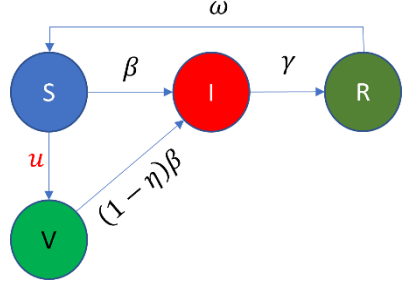

**Figure 2:** Model Flowchart (incorporating optimal control). This diagram illustrates the progression of individuals within the model. Susceptible individuals may contract the infection at a transmission rate, $\beta$, moving from compartment $S$ to compartment $I$. Infected individuals can recover at a rate of $\gamma$, transitioning to the recovered compartment $R$. Susceptible individuals may opt for vaccination, transitioning to compartment $V$ at a rate determined by optimal control, denoted by $u$. Vaccinated individuals may still become infected at a discounted transmission rate $(1 - \eta)\beta$, moving to the infected compartment. Recovered individuals may lose immunity over time, potentially becoming susceptible again at a rate of waning immunity represented by $\omega$.

**Table 1:** Description of the model parameters

| Parameter symbol | Parameter Description |
|---|---|
| $\beta$ | Disease Transmission rate |
| $\gamma$ | The recovery rate |
| $m$ | Effect of inertia when switching from $S$ to $V$ |
| $c$ | Cost of infection |
| $k$ | Relative sensitivity due to vaccine's cost |
| $c_v$ | Cost of vaccination |
| $\omega$ | Waning rate against immunity |

$$\dot{S} = -\beta SI - uS + \omega R, \tag{7}$$

$$\dot{V} = uS - (1 - \eta)\beta VI, \tag{8}$$

$$\dot{I} = \beta SI + (1 - \eta)\beta VI - \gamma I, \tag{9}$$

$$\dot{R} = \gamma I - \omega R, \tag{10}$$

$$S(t) + V(t) + I(t) + R(t) = 1, \tag{11}$$

where $u(t)$ is the vaccination control that needs to be optimized at any time $t$.

According to Pontryagin's maximum principle, we define the objective function for the above model as

$$J = \min \int_0^T (cI + c_v uS)^2 dt, \tag{12}$$

where $c$ and $c_v$ are the disease cost due to infection and the vaccination cost respectively. Note that $c \cdot I(t)$ indicates the disease cost, socially accumulated, at time $t$ and $c_v \cdot u(t) \cdot S(t)$ means the socially accumulated vaccination cost at time $t$. To make sure the following mathematical process obeys Pontryagin's maximum principle heathy, we impose a square operator to this instead of the simple accumulated cost, which is quite analogous to the concept of the Least Square Method (LSM). The square ensures the convexity of the function defined in the objective function according to Pontryagin's maximum principle.

Next, we define the Hamiltonian as follows:

$$H = (cI + c_v uS)^2 + \lambda_1 \dot{S} + \lambda_2 \dot{V} + \lambda_3 \dot{I} + \lambda_4 (\gamma I - wrR), \tag{13}$$

where

$$\dot{\lambda}_1 = -\frac{\partial H}{\partial S} = 2(cI + c_v uS)c_v u + \lambda_1(-\beta I - u) + \lambda_2 \beta I + \lambda_3 u, \tag{14}$$

$$\dot{\lambda}_2 = -\frac{\partial H}{\partial I} = 2(cI + c_v uS)c + \lambda_2(\beta S - \gamma + (1-\eta)\beta V) + \lambda_3(1-\eta)\beta V, \tag{15}$$

$$\dot{\lambda}_3 = -\frac{\partial H}{\partial V} = \lambda_2(1-\eta)\beta I - \lambda_3(-(1-\eta)\beta I), \tag{16}$$

$$\dot{\lambda}_4 = -\frac{\partial H}{\partial R} = \lambda_1 \omega - \lambda_4 \omega, \tag{17}$$

$$\lambda_i(T) = 0, \text{ (Transversality condition)} \tag{18}$$

Thus the optimality condition is,

$$\frac{\partial H}{\partial u} = 2(cI + c_v uS)c_v u - \lambda_1 S + \lambda_3 S = 0, \text{ at } u^* \tag{19}$$

which implies, $u^* = \frac{1}{c_v S}\left(\frac{\lambda_1 - \lambda_3}{2} - cI\right)$, given that $\lambda_1 - \lambda_3 \geq 2cI$, and $c_v \neq 0$. \hfill (20)

Thus, the optimal control will be

$$u^* = \min\left\{\max\left\{0, \frac{1}{c_v S}\left(\frac{\lambda_1 - \lambda_3}{2c_v} - cI\right)\right\}, u^{max}\right\}, \tag{21}$$

$u^{max}$ is the maximum rate of control that can be applied.

## 2.3 Primary Reproduction Number, Cumulative infection, Cumulative vaccination, Average Social payoff (ASP), Social efficiency deficit (SED)

In this study, we considered the primary reproduction number, $R_0 = \frac{\beta}{\gamma} = 2.5$ [4,47,48,58].

The cumulative number of infected individuals is:

$$IT = \int_0^\infty (\beta SI + (1-\eta)\beta VI)\, dt, \tag{22}$$

The cumulative number of vaccinated individuals is:

$$VT = \int_0^\infty xS\, dt, \left[\int_0^\infty uS\, dt \text{ for the optimal control}\right] \tag{23}$$

Where $t = \infty$ indicates a state of equilibrium (we say it, Nash equilibrium, NE).

The Average social payoff (ASP) at NE is described as follows.:

$$ASP_{NE} = -IT * C - VT * c_v, \quad [x \text{ is the vaccination rate}] \tag{24}$$

The first and second terms on the right show the total rewards for the infected and immunized individuals, respectively.

The ASP at Social Optimum is defined as follows:

$$ASP_{SO} = -IT * C - VT * c_v, \quad [u \text{ is the vaccination rate}] \tag{25}$$

In the model, the first expression on the right-hand side represents the sum of people's payoff who become infected. In contrast, the second term represents the total payoffs of vaccinated individuals. The SED is defined as follows:

$$SED = ASP_{SO} - ASP_{NE}. \tag{26}$$

## 3. The Findings and Discussion

### 3.1 Illustration from timeseries (with typical values):

In this section, we provide the time series data for the compartments and vaccination rates in the behavior and optimal control models using values of the common parameters. The initial values for the compartments and rates are shown in Table 3, whereas the standard parameter values are shown in Table 2. To solve equations (1)−(6) in the behavior model, we employed the explicit finite difference method with a time step size of $dt = 0.1$. We solved equations (7)−(21) for the optimal control model using the forward-backward sweep method and fourth-order Runge-Kutta method, with a time step size of 0.1. In the behavior model, we set the initial vaccination rate to 0.1, while in the optimal control model, we assumed that a maximum control rate of vaccination, $u^{max}$, could be applied up to 0.1. This value is considered reasonable and feasible by any authority.

Figure 3(a) presents the compartments' time series using the behavior model's dynamics with standard parameter values. The susceptible ($S$), infected ($I$), recovered ($R$), and vaccinated ($V$) compartments are depicted over time. Similarly, in Figure 3(b), we display the corresponding vaccination rate using the same parameter settings. Moving on to Figures 3(c) and 3(d), we showcase the time series of the model compartments and the optimal control vaccine flow, respectively, using the suggested optimal control concept. The optimal control approach notably leads to a more stabilized vaccination rate than the behavior model dynamics. Even though the susceptible population continues to increase, the vaccination rate reaches a steady 400 days. Contrastingly, Figures 3(a) and 3(b) demonstrate that without the optimal control approach, the number of vaccinated and susceptible individuals gradually increases, necessitating a prolonged vaccination campaign. Furthermore, the peak of the vaccination flow approaches 12%. Therefore, based on the same parameter combination, utilizing the optimal control idea proves more advantageous to society in regulating vaccination strategies and achieving desired outcomes.

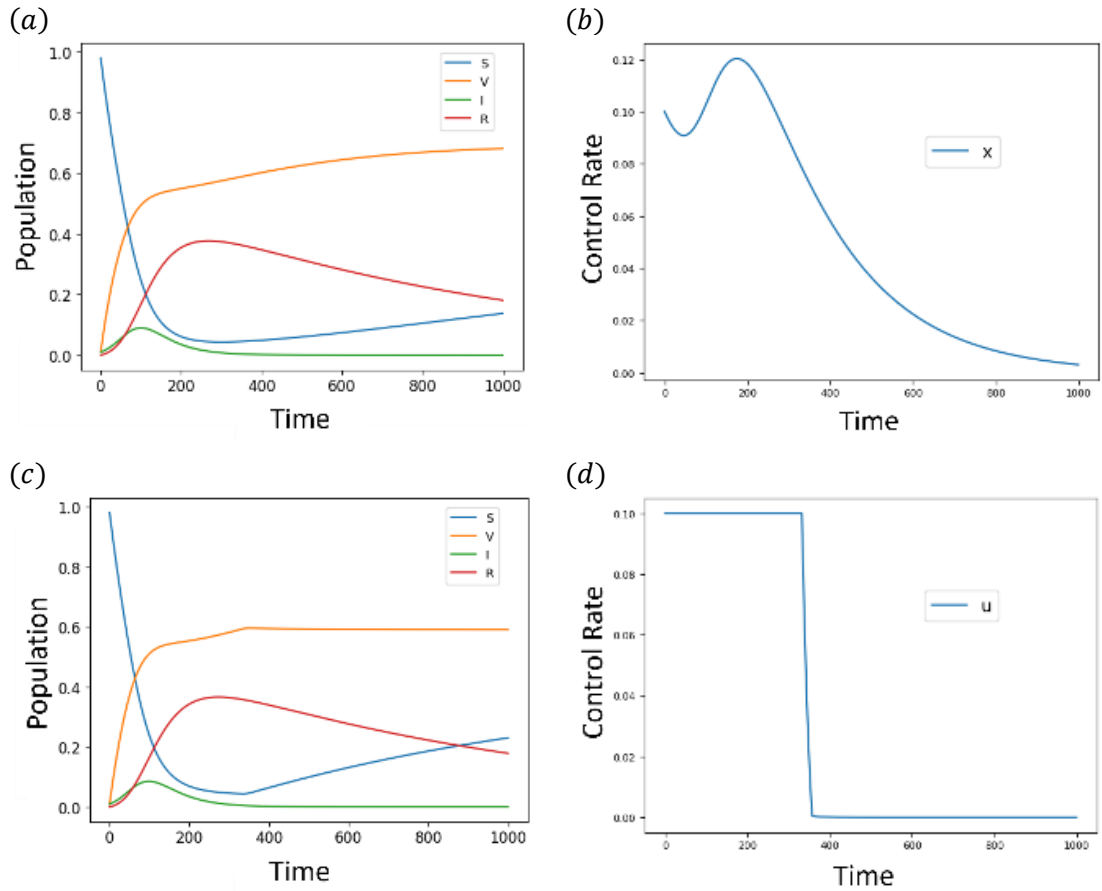

**Figure 3:** In panels (a) and (b) of the figure, the time series representation of the compartments and vaccination rate is provided for the behavior model. The blue line corresponds to the susceptible individuals, the green line represents the infected individuals, the orange line represents the immunized individuals, and the red line corresponds to the recovered individuals. In panel (b), the line depicts the vaccination rate. Similarly, panels (c) and (d) illustrate the compartments and vaccination rate for the optimal control model. The blue, green, orange, and red lines in panel (c) represent the susceptible, infected, vaccinated, and recovered individuals, respectively. Panel (d) shows the vaccination rate. By comparing the behavior model with the optimal control model, using the same parameter settings, we observe that the vaccination equilibrium reaches sooner in the optimal control model. So the optimal control strategy facilitates the achievement of a stable vaccination rate in a shorter duration compared to the behavior model.

**Table 2:** Standard values of the parameters [38,47,50,52]

| Parameter | Value | Parameter | Value |
|---|---|---|---|
| $\beta$ | 0.833 | $\omega$ | 1/90 |
| $\gamma$ | 0.333 | $k$ | 0.1 |
| $m$ | 1.0 | $c_v$ | 0.5 |
| $c$ | 1.0 | $u^{max}$ | 0.1 |

**Table 3:** Initial values for the compartments and vaccination rates

| State | At $t=0$ | State/Rate | At $t=0$ |
|---|---|---|---|
| $S$ | 0.98 | $R$ | 0.00 |
| $V, I$ | 0.01 | $x$ | 0.1 |

## 3.2 Timeseries comparison based on wanning immunity, $\omega$:

In this section, we will examine the impact of the immunity-waning rate on both the behavior model and the optimal control model. Figure 4 displays the time series of the compartments to observe the effect. We consider four different values for the waning rate of immunity ($\omega = 0.0, \frac{1}{90}, \frac{1}{60}$, and $\frac{1}{30}$ $day^{-1}$). In panels (a) to (d), the time series of the susceptible, vaccinated, infected, and vaccination rates are shown for the behavior model with varying values of the immunity-waning rate. Correspondingly, panels (e) to (h) illustrate the identical diagrams for the optimal control model. When we observe panels (a) and (e), we notice that the behavior of the susceptible cases remains almost the same across different values of $\omega$. The patterns in the susceptible compartment are similar irrespective of the waning rate of immunity. Some significant characteristics are evident in panels (b) and (f), which represent the vaccinated individuals. In the case of the behavior model, the vaccination peak is significantly higher than the optimal control model. Additionally, the optimal control model's equilibrium point occurs earlier than the behavior model's. This trend holds for all values of $\omega$. Panels (c) and (d), representing the infected individuals, exhibit essentially the same patterns in both models. Despite variations in the level of immunity, the peaks of infection remain unchanged in both the behavior and optimal control models. Finally, panels (d) and (h) depict the vaccination rates. Between the two models, noticeable differences were observed. The optimal control model ensures that vaccination is not continued until close to the end of the season, resulting in greater cost-effectiveness for the authorities and better overall outcomes. Considering all these observations, we can conclude that increasing the waning immunity rate encourages more individuals to get vaccinated, which aligns with our expectations.

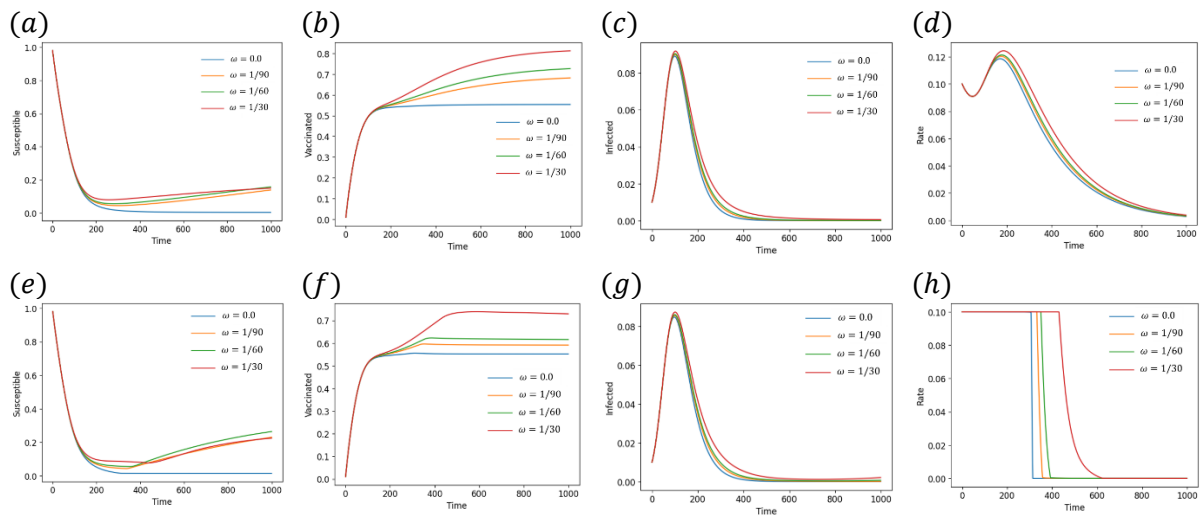

**Figure 4:** In panels (a) to (c), the time series of the susceptible, vaccinated, and infected compartments are displayed. The curves in blue, orange, green, and red correspond to the waning rates of $0.0, \frac{1}{90}, \frac{1}{60}, and \frac{1}{30}$ $day^{-1}$, respectively. These panels provide insights into the behavior model. Similarly, panels (e) to (g) show the time series of the susceptible, vaccinated, and infected compartments for the optimal control model. The same color scheme (blue, orange, green, and red) represents the waning rates. Panel (d) shows the vaccination rate for the behavior model, while panel (h) shows the vaccination rate for the optimal control model. The panels mentioned above show that increasing the waning immunity rate leads to more vaccinations. Additionally, we can see that the optimal control model achieves equilibrium more quickly than the behavior model. So the optimal control model effectively regulates the vaccination strategy and stabilizes the system at a faster pace. In summary, raising the waning rate of immunity encourages more individuals to get vaccinated, and the optimal control model outperforms the behavior model by achieving equilibrium more rapidly.

## 3.3 Time series comparison based on vaccination cost, $c_v$:

In this section, we will examine how the cost of vaccination impacts both the behavior and optimal control models. Figure 5 provides the time series of the compartments for both models to observe the effect. We consider three different values for the vaccination cost ($c_v = 0.2, 0.5,$ and $0.9$). In panels (a) to (d), the time series of the susceptible, vaccinated, infected, and vaccination rates are displayed for the behavior model with varying vaccination costs. Similarly, panels (e) to (h) illustrate the same diagrams for the optimal control model. Examining panels (a) and (e) for the susceptible cases, we observe that in the behavior model, higher vaccination costs result in a larger portion of the population remaining susceptible for a longer period. In contrast, the optimal control model shows minimal impact on the susceptible individuals as vaccination costs increase. We observe similar trends in panels (b) and (f), which represent the vaccinated individuals. Lower vaccination costs lead to higher vaccination rates in both the behavior and optimal control models. However, the optimal control model achieves equilibrium sooner than the behavior model. Panels (c) and (g) depict the infection level. In the behavior model, when the cost of vaccination increases, it leads to a rise in the infection peak. On the other hand, the optimal control model maintains a relatively steady level of infection despite variations in vaccination costs. Finally, panels (d) and (h) display the vaccination rates. Significantly different patterns emerge between the behavior model and the optimal control model. In the optimal control model, an increase in vaccination cost allows for a shorter duration of control, which is more favorable for the healthcare authority. Conversely, the behavior model exhibits a significant increase in vaccination with a decrease in vaccination cost, but this process continues until the end of the season, which can burden the healthcare authority. Considering all aspects, we can conclude that the optimal control model is more suitable for limiting the threat of the epidemic. Higher vaccination costs tend to decrease the vaccination rate, while lower costs increase the vaccination rate.

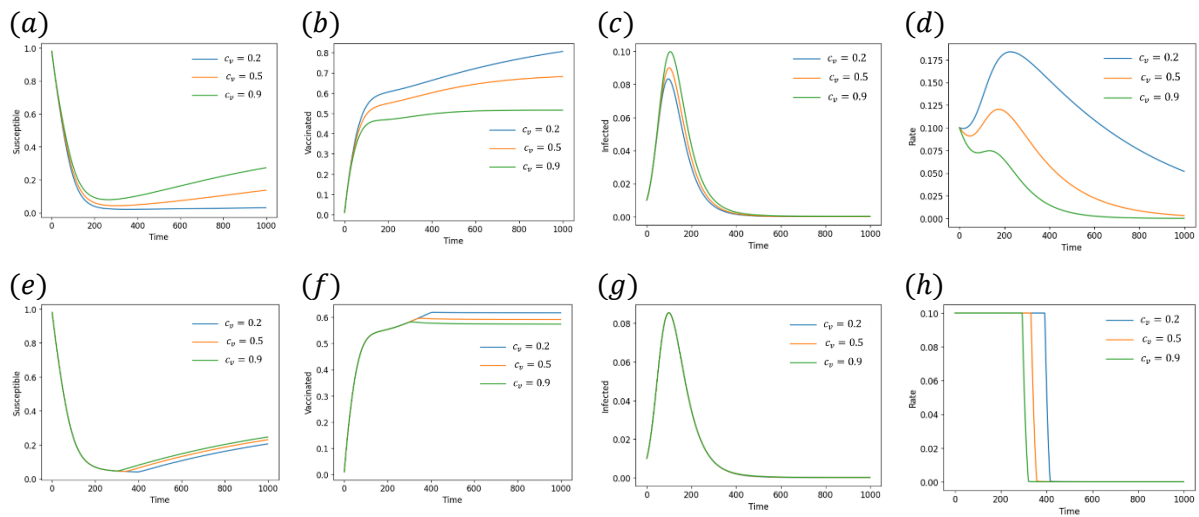

**Figure 5:** The time series of the susceptible, vaccinated, and infected compartments are depicted in panels (a) to (c) in the behavior model, respectively. Panel (d) displays the vaccination rate in the behavior models, with the blue, orange, and green curves representing vaccine costs of 0.2, 0.5, and 0.9, respectively. In panels (e) to (g), the time series of the susceptible, vaccinated, and infected compartments are presented for the optimal control model, while panel (h) shows the vaccination rate. We use the same color scheme to represent the vaccine costs. Considering the overall situation, we can conclude that optimal control is more suitable for reducing the threat of the epidemic. Higher vaccine costs lead to reduced vaccination rates, while lower costs result in increased vaccination rates. This conclusion is highly plausible and aligns with our understanding.

## 3.4 Timeseries comparison based on vaccine efficiency, $\eta$:

In Figure 6, we examine how the effectiveness of vaccination impacts the models. The time series of the compartments for the behavior model and the optimal control model is displayed to observe the effect. We consider three vaccination efficiency values: $\eta = 0.4, 0.7,$ and $0.9$. Panels (a) to (d) illustrate the time series of the susceptible, vaccinated, infected, and vaccination rates, respectively, with variations in vaccine effectiveness and the behavior model. Panels (f) to (h) present the same diagrams for the optimal control problem. Comparing panels (a) and (e) for the susceptible cases, we observe similar trends in both models, indicating a decrease in the number of susceptible individuals over time as vaccine effectiveness increases. The panels for vaccinated individuals (b) and (f) also exhibit similar patterns. Lower vaccine effectiveness leads to lower vaccination rates, but the optimal control model reaches equilibrium faster than the behavior model. The infection peak patterns are identical for both the behavior and optimal control models, as shown in panels (c) and (g). The vaccination rate panels (d) and (h) demonstrate noteworthy differences. With improved vaccine effectiveness, the optimal control model can apply control for a shorter duration, which is generally advantageous for the healthcare authority. In contrast, according to the behavior model, vaccination rates increase significantly as vaccine efficacy declines, and this process continues until the end of the season, placing a burden on the healthcare authority. Considering the complete picture, we can conclude that optimal control is more effective in reducing the threat, and higher vaccine efficiency results in fewer overall vaccinations, while lower efficiency leads to more overall vaccinations.

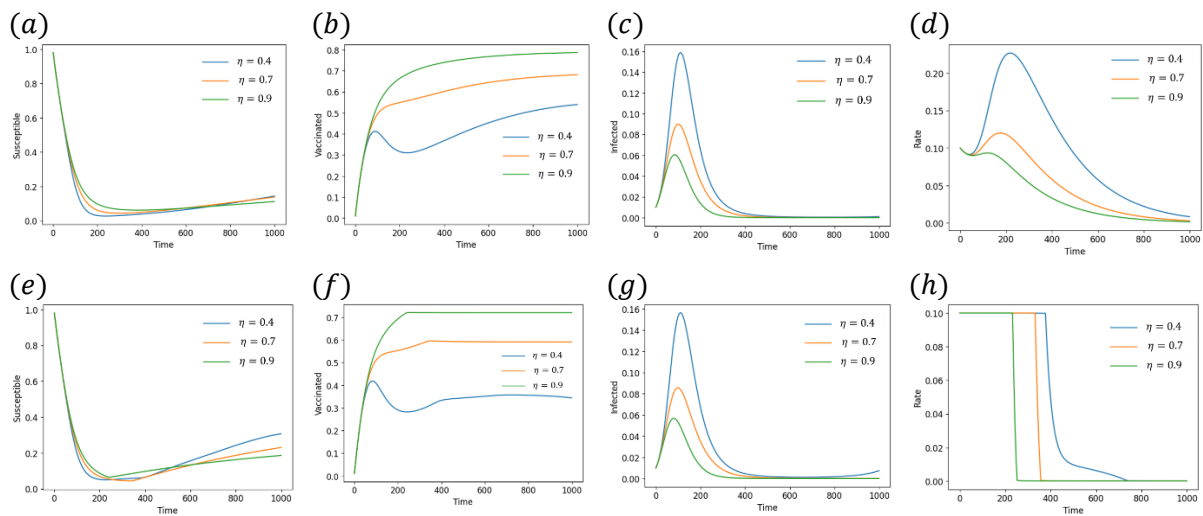

**Figure 6:** In panels (a)−(c) of the behavior model, the blue, orange, and green curves represent the time series of the susceptible, vaccinated, and infected compartments, respectively, for different values of vaccine efficacy ($\eta = 0.4, 0.7,$ and $0.9$). We can observe that as the vaccine efficacy increases, the number of susceptible individuals decreases over time. Similarly, in panels (e)−(g) of the optimal control model, the blue, orange, and green curves represent the time series of the susceptible, vaccinated, and infected compartments, respectively, with different values of vaccine efficacy. The trends are consistent with the behavior model, indicating that higher vaccine efficacy leads to fewer susceptible individuals. Panel (d) of the behavior model shows the rate of vaccination. We can see that as the vaccine efficacy increases, the vaccination rate decreases. This suggests that higher vaccine efficacy leads to a decreased requirement for vaccination. Panel (h) of the optimal control model displays the vaccination rate. Similar to the behavior model, the vaccination rate decreases as the vaccine efficacy increases. This implies that higher vaccine efficacy leads to a decrease in the optimal control's recommendation for vaccination. Considering the overall situation, We can conclude that optimal control is more appropriate for mitigating the threat of the epidemic. Additionally, higher vaccine efficacy is associated with a decrease in overall vaccination rates, while lower efficacy

promotes vaccination. This conclusion aligns with the expected behavior, as higher-efficacy vaccines would provide better protection and reduce the need for vaccination.

### 3.5 ASP and SED:

In Figure 7.1, the first row of panels (a) through (c) represents the total number of infected, vaccinated, and asymptomatic individuals at the equilibrium state using the behavior model (Nash equilibrium, NE). The second row of panels (d) through (f) represents the same quantities using the optimal control problem (Social optimum, SO). Panel (g) depicts the SED, which distinguishes the two total social payoffs in panels (f) and (c). In panel (a), which corresponds to the behavior model, the total number of infected individuals shows a monotonic increase as the transmission rate ($\beta$) increases. This indicates that higher transmission rates lead to higher infection rates in the behavior model.

Similarly, panel (d) shows that high infection rates can result from high transmission rates and low vaccine efficiency in the optimal control model. Panel (b) demonstrates that the behavior model consistently recommends vaccination throughout the season, regardless of the transmission rate. On the other hand, panel (e) shows that the optimal control model advises more vaccination as the transmission rate increases. This suggests that the optimal control model acknowledges the necessity of increasing vaccination rates in response to higher transmission rates. In panels (c) and (f), the ASPs are depicted for the behavior and optimal control models, respectively. The behavior model maintains a high level of social payoff across the entire parameter space, indicating that it is less sensitive to changes in vaccine efficiency compared to the optimal control model. The optimal control model, on the other hand, achieves higher social payoffs by balancing the trade-off between vaccination and infection. Panel (g) displays the SED, representing the difference between the optimal control model's social payoff and the behavior model's. Three regions are observed in the SED panel. The lower blue zone represents a low dilemma, where the transmission rate remains low, and both models achieve relatively high social payoffs. The dark blue area in the middle indicates almost no dilemma, as the optimal control model maintains higher vaccination rates and lower infection rates than the behavior model. However, as the transmission rate and vaccine efficiency increase, the light blue zone emerges, representing the most challenging region. In this region, the behavior model's social payoff is significantly higher than that of the optimal control model, as the behavior model is less sensitive to vaccine efficiency.

In conclusion, the optimal control model demonstrates a better balance between vaccination and infection rates, leading to higher social payoffs. While consistently recommending vaccination, the behavior model may result in higher infection costs compared to the optimal control model. The SED analysis highlights the regions where challenges and dilemmas arise, with the behavior model showing higher costs and the optimal control model achieving a better balance.

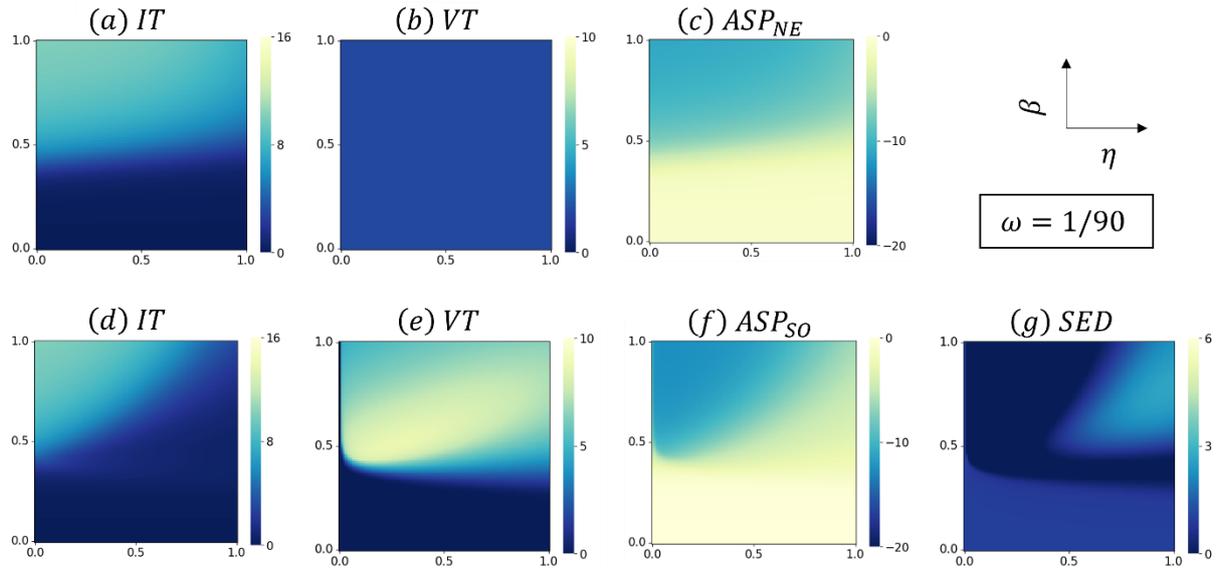

**Figure 7.1:** In the given Figure, panels (a) through (c) represent the overall prevalence of infection, vaccination, and ASP (Average Social Payoff) in the Nash equilibrium (NE) situation using the behavior model. Panels (d) through (f) represent the same quantities in the social optimum (SO) situation using the optimal control model. Panel (g) displays the SED, which distinguishes the two ASPs in panels (f) and (c). The y-axis represents the transmission rate, ranging from 0.0 to 1.0, and the x-axis represents the vaccine effectiveness, from 0.0 to 1.0. The color scale indicates the values of the respective variables, with the cumulative number of infected individuals ranging from 0 to 16, the total number of vaccinated individuals ranging from 0 to 10, and the SED ranging from 0 to 6. The ASPs ranges are $-20$ to $0$. By examining the figure, we can observe that as both the transmission rate and vaccine effectiveness increase, the social problem becomes more significant. This implies that higher transmission rates and more effective vaccines lead to more significant challenges in achieving a desirable social outcome. The SED panel (g) visually represents the differences between the optimal control model and the behavior model in terms of social payoffs. The figure demonstrates that addressing the social problem becomes more challenging as the transmission rate and vaccine effectiveness increase. This emphasizes the significance of identifying optimal control strategies to effectively manage vaccination efforts and prevent infections, ultimately leading to the attainment of optimal social outcomes.

In Figures 7.2 and 7.3, the panels from Figure 7.1 are replicated with changes in the waning immunity rate. Figure 7.2 uses a waning immunity rate of $\omega = \frac{1}{60} \, day^{-1}$, while Figure 7.3 uses a waning immunity rate of $\omega = \frac{1}{30} \, day^{-1}$. The purpose is to examine the impact of waning immunity on the social dilemma scenario while maintaining an identical range for the panels. Upon comparing the relevant panels in Figures 7.1, 7.2, and 7.3, we can observe that more individuals are vaccinated as the waning immunity rate increases. This decreases the average social payoff towards the socially ideal level in the context of optimal control. The SED panels allow us to identify regions similar to those in Figure 7.1 but with higher positive SED values. From these observations, we can conclude that increasing the rate of waning immunity amplifies the social dilemma, which aligns with expectations. This indicates that when waning immunity occurs faster, balancing vaccination efforts and achieving the desired social outcomes becomes more challenging.

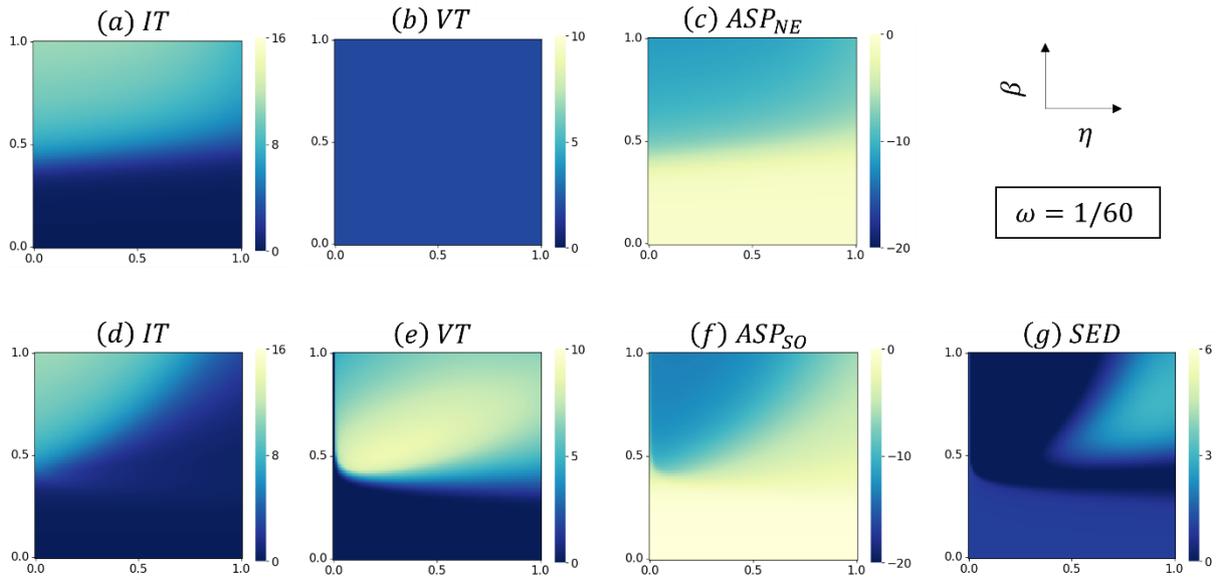

**Figure 7.2:** The Figure shows the overall prevalence of infection, vaccination, and ASP from the behavior model and optimal control model, as well as SED. All the ranges are the same as in Figure 7.1, with standard values of parameters as well, except the wanning immunity rate ($\omega$).

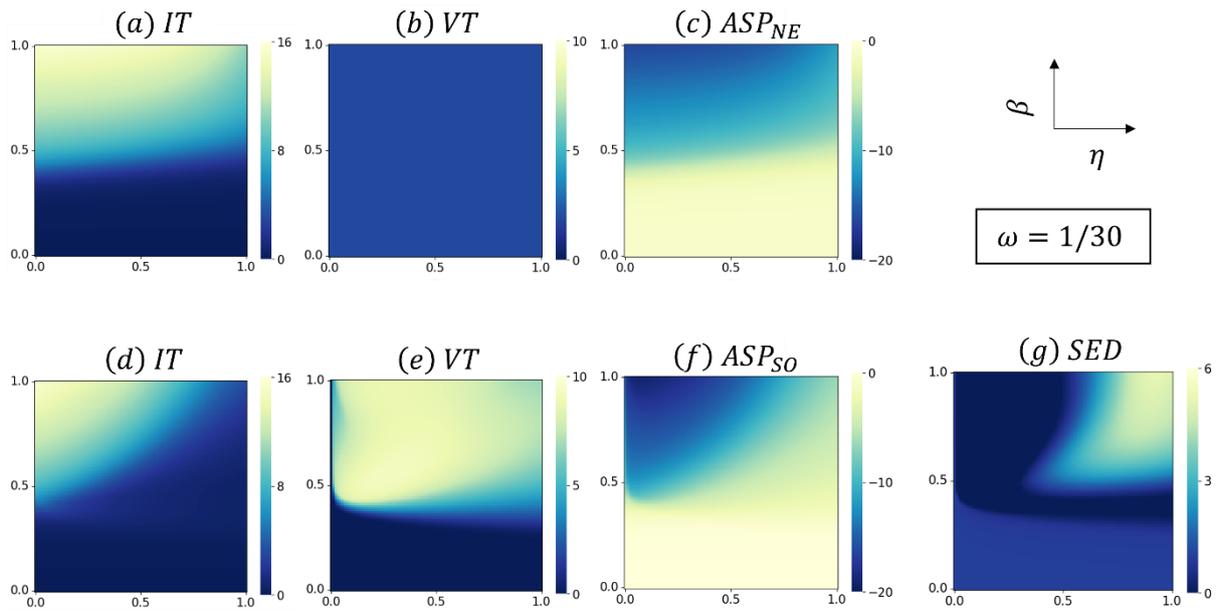

**Figure 7.3:** This Figure displays the overall prevalence of infection, vaccination, and ASP from the behavior model and optimal control model, as well as SED. Except for the wanning immunity rate ($\omega$), all ranges and standard parameter values are identical to Figure 7.1.

Figure 7.4 focuses on the impact of the cost of the disease ($c$) and the cost of vaccination ($c_v$) on the social dilemma scenario, as well as the total number of infected individuals, vaccinated individuals, and the ASPs for both the behavior and optimal control models. The figure considers a range of values for the disease cost ($c$) from 0.01 to 1.01 and the vaccination cost ($c_v$) from 0.01 to $c$, with all other parameters set to their default values. From the panels representing the behavior model in the first row, we can observe that as the cost of vaccination increases, fewer individuals choose to get vaccinated, resulting in a higher overall number of people becoming infected. However, when comparing the behavior model to the optimal control model, it is evident that vaccination rates remain high and regional infection rates remain low in the latter. The panels depicting the average social payoff show a similar trend. When the costs of disease and vaccination are combined, the average social payoff increases as the cost of vaccination rises. This suggests higher vaccination costs incentivize individuals to prioritize vaccination, improving overall social outcomes. One interesting observation can be made from the SED panel. The SED is largest when the disease cost is high, and the cost of vaccination is around 0.4. This indicates that when the cost of the disease is high, but the cost of vaccination is moderately high, people hesitate to get vaccinated. As a result, the decision to get vaccinated becomes more dependent on the expense of the vaccination itself. Similar patterns have been observed in previous figures, such as Figures 7.1-7.3, where increasing the immunity rate resulted in higher values of SED, indicating an intensified social dilemma. This behavior aligns with human behavior and decision-making processes.

In summary, Figure 7.4 demonstrates how the costs of disease and vaccination impact vaccination rates, infection rates, average social payoffs, and the intensity of the social dilemma. Higher vaccination costs lead to lower and higher infection rates, while the combined costs of disease and vaccination influence the average social payoff. The SED panel reveals the regions where the social dilemma is most pronounced, highlighting the interplay between vaccination costs and disease in shaping individual behavior and overall social outcomes.

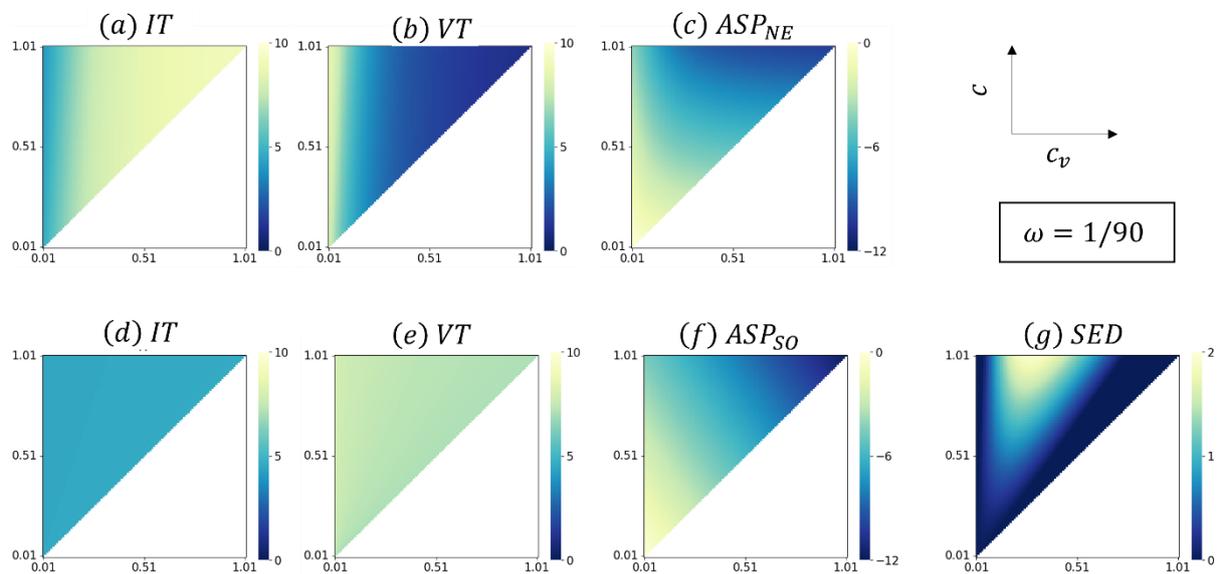

**Figure 7.4**: This Figure illustrates the prevalence of infection, vaccination, and the ASPs from both the behavior model (NE) and the optimal control model (SO), along with the SED panel. Panels (a) through (c) represent the NE condition, while panels (d) through (f) depict the SO case. Panel (g) specifically shows the SED, distinguishing it from panels (f) and (c). The figure uses the cost of disease ($c$) on the y-axis and the cost of vaccination ($c_v$) on the x-axis, with ranges of 0.01 to 1.01 for the disease cost and 0.01 to $c$ for the vaccination cost. The cumulative number of infected and vaccinated people ranges from 0 to 10, while the ASP ranges from -12 to 0. The SED ranges from 0 to 2, and the remaining variables follow the base case. As the cost of disease increases, it is evident that the social dilemma worsens. However, the most significant social dilemma occurs when the cost of vaccination is approximately 0.4, and it diminishes as the vaccination cost increases. This implies that when the cost of vaccination is moderate, people are more hesitant to get vaccinated despite the higher cost of the disease. This behavior highlights the delicate balance between vaccination costs and disease in influencing individual decisions and social outcomes.

In summary, Figure 7.4 demonstrates the relationship between disease costs and vaccination and their impact on the social dilemma. Higher costs of disease exacerbate the social dilemma, but the intensity of the dilemma is greatest when the cost of vaccination is around 0.4. As the cost of vaccination increases, the social dilemma becomes less pronounced. The SED panel visually represents the regions where the social dilemma is most significant, shedding light on the interplay between the costs of vaccination and disease in shaping individual choices and the overall social dynamics.

## 4. Conclusion

When exploring an epidemic model and its associated social dilemma, particularly concerning vaccination, the primary focus is on implementing effective vaccination strategies. For individuals, the priority lies in vaccination uptake, taking into account factors like the eventual epidemic size, its peak, and the cost of vaccination. Conversely, authorities aim to minimize disease spread while keeping costs low. The rate of immunity waning significantly affects strategy implementation for both individuals and authorities. Thus, integrating the epidemic model with optimal control theory becomes a valuable tool for healthcare professionals and management authorities in crafting and sustaining effective vaccination strategies. Previous studies have offered various analytical and numerical insights into optimal vaccination control strategies, often utilizing Pontryagin's maximum principle [7,8,11,15,23,24,26,37,56]. However, these studies frequently overlooked vaccination costs in their models and analyses. Moreover, some explored constant-rate vaccination or interventions, posing challenges for real-world implementation [4,38,47,52,58,59,62].

This research takes a direct approach by incorporating vaccination costs into the analysis to derive an optimal vaccination control strategy while upholding Pontryagin's maximum principle. Unlike prior works, we introduce a novel objective function that rigorously reflects the socially accumulated total cost. We compare a conventional cost-based behavior model, which considers human responses to vaccination based on epidemic conditions and costs, with our proposed optimal control model integrating vaccination and infection costs under the same parameters. Both models undergo thorough analysis and comparison, revealing the practicality and quicker stabilization of vaccination with the optimal control model. Regarding the social dilemma, numerical results illustrate how increasing rates of immunity waning amplify the Social Efficiency Deficit (SED), while higher vaccination costs somewhat mitigate it. Furthermore, the method used to calculate social dilemmas in this research is deemed more reliable for understanding social scenarios.

While this study focuses on a simple model incorporating a single control—vaccination—it's crucial to note that managing pandemics involves considering various interventions like isolation, quarantine, treatment, and testing policies. Future work will extend our model to include multiple interventions to enhance the effectiveness of analyzing epidemic models and optimal control strategies for preventive measures. This broader perspective aims to provide a comprehensive understanding of how different interventions interact to manage and control pandemics.

**Acknowledgment**

Grant-in-Aid for Scientific Research from JSPS, Japan, KAKENHI (Grant No. JP 23H03499), given to Professor Tanimoto, provided a portion of the funding for this study. We wish to thank them for everything.

**Data availability**

No data was used for the research described in the article.

**Declaration of Conflict of Interest**

No conflict of interest exists.

## Authors Contribution

**Md. Mamun-Ur-Rashid Khan**: Conceptualization, Formal analysis, Investigation, Methodology, Software, Validation, Visualization, Writing – original draft, Writing – review & editing.

**Jun Tanimoto**: Conceptualization, Formal analysis, Funding acquisition, Supervision, Writing – original draft, Writing – review & editing.